\documentclass[twoside,12pt]{article}
\usepackage{amsthm,amsmath,amsfonts,amscd,amssymb,dsfont}
\usepackage{latexsym,amsmath}
\usepackage[T1]{fontenc}
\usepackage[pagebackref=false,colorlinks,linkcolor=black,citecolor=magenta]{hyperref}
\usepackage{framed}
\usepackage{pgfplots}
\usepackage{tikz-cd}
\usepackage{mathrsfs}
\usepackage[all]{xy}
\usetikzlibrary{arrows}
\usepackage{hyperref}
\usepackage{mathrsfs}
\usepackage{upgreek}
\usepackage{blkarray}
\usetikzlibrary{arrows}
\usepackage{pgfplots}
\pgfplotsset{compat=1.7}
\usepackage{etoolbox,lastpage}
\usepackage{comment}
\usepackage{graphicx}
\usepackage{enumitem}
\usepackage[titletoc,toc]{appendix}
\newtheorem{theorem}{Theorem}[section]
\newtheorem{lemma}[theorem]{Lemma}
\newtheorem{proposition}[theorem]{Proposition}

\newtheorem{definition}[theorem]{Definition}
\newtheorem{remark}[theorem]{Remark}
\newtheorem{example}[theorem]{Example}
\numberwithin{equation}{section}
\allowdisplaybreaks

\makeatother
\begin{document}
	\title{\large\bf
		Euler-Poincare characteristic pair of orientable supermanifolds}
	\bigskip
	\author{{\small \bf Mehdi Ghorbani $^1$, Fatemeh Alikhani $^1$, Saad Varsaie $^1$\footnote{
				{ Corresponding author, E-mail address: \url{varsaie@iasbs.ac.ir}}\\$^1$ Institute for Advanced Studies in Basic Sciences, 444 prof. Yousef Sobouti Blvd, Zanjan 45137-66731, Iran\\}}~\\[5mm]}
	\date{}
	\maketitle	
	\begin{abstract}
		The Euler-Poincare characteristic, or Euler characteristic in short, is a fundamental topological invariant of compact manifolds that plays a crucial role in a variety of geometric and topological situations. From this point of view, we tried to expand on this important concept in supergeometry. In this article, we introduce the Euler-Poincare characteristic pair in the supergeometry. In the final section, we examine transversality in the category of
		$\Pi$-symmetric supermanifolds. 
	\end{abstract}
	\bigskip
	\baselineskip=0.2cm
	\tableofcontents
	\bigskip
	\baselineskip=0.5cm
	\section{Introduction}
	Euler characteristic is a basic invariant for compact manifolds. It is defined not only for manifolds but also for real vector bundles whose rank equals the dimension of their base manifold. It plays a crucial role in both algebraic and differential topology. It is very difficult to calculate the Euler characteristic using cohomology in algebraic topology. However, there is a simple way to calculate the Euler number in differential topology, which is expressed using a self-intersection number.
	In Morse theory perspective, the Euler characteristic can be interpreted as The Poincaré-Hopf theorem relates to the Euler characteristic of a manifold and critical points of a smooth function. In differential topology relates Euler characteristic with the transversal intersection number between a global section and zero section of a vector bundle on a smooth manifold.
	In fact, the Euler characteristic is a bridge between algebraic topology and differential topology.\\
	The Poincaré-Hopf theorem states that for a smooth vector field on a closed manifold M, the sum of the indices of its zeros equals the Euler characteristic of M. The Euler characteristic in Gauss-Bonnet Theorem For closed Riemannian manifolds can be calculated by integrating the curvature. In algebraic geometry, the Euler characteristic of a coherent sheaf F on a proper scheme X is defined as the alternating sum of the dimensions of its cohomology groups. In differential geometry, the Euler characteristic is typically defined as the alternating sum of the Betti numbers of a topological space.\\ 
	Extending our prior work\cite{4} on transversality and genericity within the category of supermanifolds, the present article develops a rigorous framework for defining and computing self-intersection numbers in this context. Then we derive a natural formulation of the Euler–Poincare characteristic as a homotopy-invariant intersection number. This construction not only generalizes classical results from differential topology to the supergeometric setting, but also establishes foundational tools for further advances in index theory, topological invariants in super field theories, and the intersection theory of derived moduli spaces.
	\section{Preliminaries}
	\textbf{Supermanifold}.
	A supermanifold of dimension $m|n$ is a ringed space 
	$X=(\tilde{X},\mathcal{O}_X)$ that is locally isomorphic to $U^{m|n}$ and $\tilde{X}$ is a second countable
	and Hausdorff topological space. The superdomains $\mathbb{R}^{m|n}$ and $U^{m|n}$ are special examples of supermanifolds.\\
	\textbf{Morphism}.
	A morphism $\varphi$ from the supermanifold $X=(\tilde X,\mathcal O_X)$ to the supermanifold $Y=(\tilde Y,\mathcal O_Y)$ consists of a pair $(\tilde \varphi, \varphi^*)$, where $\tilde \varphi:\tilde X\to \tilde Y$ is a continues function and
	$\varphi^*: \mathcal O_Y \to \varphi_*(\mathcal O_X)$
	is a homomorphism between sheaves of supercommutative $\mathbb Z_2$- graded rings. In other words, for each open  $\tilde V \subset \tilde Y$, there is a homomorphism
	$$\varphi^*_{\tilde V}: \mathcal O_Y (\tilde V) \to \mathcal O_{X} ( \tilde \varphi^{-1} (\tilde V))$$
	such that for every open subsets $\tilde U$ and $\tilde V$ with $\tilde U\subset \tilde V$ the following diagram commutes
	\begin{equation}
		\begin{CD}
			\mathcal{O}_Y(\tilde V)@>\varphi^*_{\tilde V}>>\mathcal{O}_X(\tilde \varphi^{-1} (\tilde V))@.\qquad \\
			@ Vr^Y_{\tilde U\tilde V}VV @ VV r^X_{\tilde U\tilde V} V @.\\
			\mathcal{O}_Y(\tilde U)@>>\varphi^*_{\tilde U}>\mathcal{O}_X(\tilde \varphi^{-1} (\tilde U))@.
		\end{CD}.
	\end{equation}
	This collection of homomorphisms can be thought of as a set of constraints that relate the structure sheaf of $X$ to the structure sheaf of $Y$. It is worth noting that by these constraints it can be shown that $\tilde \varphi$ is a smooth chart if $\tilde{X}$ and $\tilde{Y}$ are equipped with their differential structures. Indeed, One can show that there exists a unique differential structure on $\tilde X$ under which $\tilde X$ is a smooth manifold embedded in $X$.\\
	\textbf{Gluing}.\label{gluing}
	Let $X_i$ be an open cover of $X$ and let $\mathcal{O}_i$ be a sheaf of rings on $\tilde X_i$. In addition assume $ f_{ij}: (X_{ji},\mathcal{O}_j|_{X{ji}})\to (X_{ij},\mathcal{O}_i|_{X{ij}})$ is an isomorphism of ringed spaces with $\tilde{ f_{ij}}=id|_{X_{ij}}$, where $X_{ij}=X_i \cap X_j$ for each $i,j.$ We call a ringed space $(X, \mathcal{O})$ is constructed by gluing ringed spaces $(X_i, \mathcal{O}_i)$ through $f_{ij}$ if there exists ringed space isomorphisms $$f_i: (X_i, \mathcal{O}|_{X_i})\to(X_i, \mathcal{O}_i|_{X_i})$$ with $\tilde{f}_i=id_{X_i}$ such that $f_{ij}=f_i\circ f_j^{-1}$.\\ 
	\textbf{Embedding}. 
	Let $\iota:=(\tilde{\iota},\iota^*):X=(\tilde{X},\mathcal{O}_X) \to Y=(\tilde{Y},\mathcal{O}_Y) $ be a morphism of supermanifolds. We say that  $(X,\iota)$ is an embedded supermanifold if $\tilde{\iota}$ is an immersion and $\tilde{\iota}:\tilde{X}  \to \tilde{Y} $ is an homeomorphism onto its image.\\
	In particular, if $\iota  (X)\subset Y $ is a closed subset of $Y$ we will say that $(X,\iota)$ is a closed embedded supermanifold.\\
	In what follows, we will always deal with closed-embedded supermanifolds. Remarkably, it is possible to show that a morphism  $\iota:X=(\tilde{X},\mathcal{O}_X) \to Y=(\tilde{Y},\mathcal{O}_Y)$ is an embedding if and only if the corresponding morphism $\iota^*:\mathcal{O}_Y \to \mathcal{O}_X $ is a surjective morphism of sheaves. Notice that, for example, given a supermanifold $Y=(\tilde{Y},\mathcal{O}_Y) $, one always has a natural closed embedding: the map $\iota:(\tilde{Y},\mathcal{O}_Y)_{red} \to (\tilde{Y},\mathcal{O}_Y) $, that embeds the reduced manifold underlying the supermanifold into the supermanifold itself.\\
	\textbf{Submersion} \label{2.9.}
	A morphism $\psi$ is called a submersion at $x$ if $d\psi_x$ is surjective. A standard submersion is a morphism $\psi: V^{m+r|n+s} \to U^{m|n}$, with coordinates $(x^i,y^a,\xi^j,\theta^b) $ and  $(x^i,\xi^j)$ for $V$ and $U$, respectively such that $\tilde{\psi}: (x, v)\mapsto x $ and
	$$\psi^*:x^i \mapsto x^i, \, \xi^j \mapsto  \xi^j. $$
	For more details on the geometry of supermanifolds, the reader is suggested to
	refer to the book \cite{3,6}.
	\section{Orientability of supermanifolds}
	In this section, the concept of orientability in supermanifolds is examined in detail from the reference \cite{1}
	and equivalent propositions are presented for it.
	\begin{definition}
		Let
		$U^{m|n}$
		be a superdomain. If
		$(x_1,...,x_m;\xi_1,...,\xi_n)$
		and
		$(y_1,...,y_m;\eta_1,...,\eta_n)$ be
		two coordinate systems on
		$U$, if for  each
		$ x \in \tilde {U} $,
		$$det\left( {\dfrac{\partial x}{\partial y}}\right)_x^\sim >0\,\,\,\text{and}\,\,\,det\left( {\dfrac{\partial \xi}{\partial \eta}}\right)_x^\sim >0$$
		then, the two coordinate systems have a similar orientation.
		The orientation
		$\epsilon$
		on this superdomain is a class of coordinate systems with similar orientations. Each superdomain U with an orientation $\epsilon$
		is denoted by 
		$(U,\epsilon)$.
	\end{definition}
	\begin{remark}
		Two oriented charts
		$(U,\epsilon_U)$ and
		$(V,\epsilon_V)$
		are intersected if
		$U \cap V$
		is nonempty and orientations coincide at the intersection.
	\end{remark}
	\begin{definition}
		Suppose
		$M^{m|n}$
		be a supermanifold. By gluing charts
		$(U,\epsilon)$ orienting covering
		$P$
		for $M$ is constructed. If $P$ has $4$ components then $M$ is orienatable\cite{1}.
	\end{definition}
	In the next section, we have introduced $P$ in more detail.
	\subsection{Orienting covering space}\label{1.1}
	Let
	$M=(\tilde M , \mathcal O_M)$
	be a supermanifold with an oriented body. Let
	$\{U_i, \epsilon_i\}_i$
	be a set of all oriented charts of
	$ M$.
	We construct orienting covering space
	$(P, \pi,M)$
	by gluing
	$(U_i, \epsilon_i)$.
	Since
	$\tilde M$
	is oriented, so according to the classical case, its orienting covering space $\tilde{P}$ is a disconnected space with two components $P_0$
	and $P_1$. Let us consider the associated
	covering map
	$$\tilde \pi: \tilde P \to \tilde M.$$
	For a local coordinate chart $U_i$ one has:
	$$\tilde \pi^{-1}\tilde{(U_i)}=\tilde V_{i0}\cup \tilde V_{i1}$$
	and
	$\pi|_{V_{ir}}$,
	for
	$r=0,1$
	is a homeomorphism.
	Now let
	$(x_1,...,x_m;\xi_1,...,\xi_n)$
	be a local coordinate system on $U_i$ that specifies orientation class
	$\epsilon_i$.
	Then the following oriented charts are located in
	$P$:
	\begin{align*}
		&(U_i, \epsilon_{i(0, 0)}),
		\,\,\,(U_i, \epsilon_{i(0, 1)}),
		\,\,\,(U_i, \epsilon_{i(1, 0)}),
		\,\,\,(U_i, \epsilon_{i(1, 1)})
	\end{align*}
	where
	$\epsilon_{i(a, b)}$
	is orientation induced by coordinates
	$((-1)^ax_1,x_2 ,...,x_m,;(-1)^b\xi_1,\xi_2,...,\xi_n)$.
	Suppose oriented charts $V_{i(a, b)}$ and $V_{j(c, d)}$ are intersected.
	Let $V_{i(a, b)}= (U_i, \tilde{\epsilon}_{i(a, b)})$
	where $\tilde{\epsilon}_{i(a, b)}$ is the orientation class
	specified by 
	$$((-1)^ax_1, 
	..., x_m; 
	(-1)^b\xi_1, ..., \xi_n).$$
	Suppose oriented charts
	$ V_{i(a, b)}$
	and
	$ V_{j(c, d)}$
	are intersected. If $f_i$ and $f_j$ are isomorphisms corresponding coordinate systems
	$$((-1)^ax_1,...,x_m;(-1)^b\xi_1,...,\xi_n)$$
	on
	$\mathcal{O}_{M}(U_i)$
	and
	$$((-1)^cy_1,...,y_m;(-1)^d\eta_1,...,\eta_n)$$
	on
	$ \mathcal{O}_{M}(U_j)$ respectively, then $f_{ij}=f_j\circ f_i^{-1}$, can be considered as transition morphism.
	Obviously,
	$f_{ij}$
	satisfy gluing conditions so oriented charts glue to each other and construct the supermanifold
	$P$.\\
	If the supermanifold
	$P= (\tilde P,\mathcal O_P)$ has 4 connected components, then $M$ is called oriented. \\
	In the sequel, the concept of the connected component for the supermanifold $P$ is explained. The connected component $P$ is a set containing oriented charts that are either pairwise intersected or there exists a finite sequence $\{V_i\}$ of oriented charts between them such that $V_i$ and $V_{i+1}$ are pairwise intersected.\\
	For the denomination of these components, one needs to define the group action of  
	$\mathbb Z_2 \oplus \mathbb Z_2$
	on oriented charts as follows:
	\begin{align}\label{1}
		T:\mathbb Z_2 \oplus \mathbb Z_2 .P \nonumber
		&\to P\\
		(a,b).(U_i,\epsilon_i)&\mapsto (U_i,\epsilon_{i(a,b)})
	\end{align}
	where $\epsilon'$ is an orientation class specified by 
	$$(a,b).(x_1,...,x_m;\xi_1,...,\xi_n):=\left( (-1)^ax_1,x_2, ...,x_m;(-1) ^b\xi_1, \xi_2, ...,\xi_n\right)$$
	if $(x_1, ..., x_m; \xi_1, ..., \xi_n)$ is a representative of orientation class $\epsilon$. In other words, here it is assumed that coordinate vector fields corresponding to the coordinate system
	$(x_i; \xi_j)$
	is a representative for orientation class $\epsilon$ and coordinate vector fields corresponding to the following coordinate system
	$$((-1)^a x_1, x_2, ...x_m; (-1)^b\xi_1, \xi_2, ...,\xi_n)$$
	is a representative for the orientation class
	$\epsilon'$. 
	Now for a specific oriented chart
	$ (U, \epsilon)$, let denote the connected component containing oriented chart
	$(a, b).(U, \epsilon) $
	with
	$P_a^b$. According to the definition, if all these components are different from each other, then
	$M$
	is called an oriented supermanifold. Assuming orientability of
	$M$,
	the group action of
	$\mathbb Z_2\oplus \mathbb Z_2$
	on these components can be expressed by the following table:
	\begin{center}
		\begin{tabular}{|c|c|c|c|c|}
			\hline
			\text{act}& $P_1^1$ & $P_1^0$ & $P_0^1$ & $P_0^0$\\
			\hline
			\text{(0,0)} & $P_1^1$ & $P_1^0$ & $P_0^1$ & $P_0^0$  \\
			\hline
			(0,1) & $P_1^0$ & $P_1^1$ & $P_0^0$ & $P_0^1$  \\
			\hline
			(1,0) & $P_0^1$ & $P_0^0$ & $P_1^1$ & $P_1^0$ \\
			\hline
			(1,1) & $P_0^0$ & $P_0^1$ & $P_1^0$ & $P_1^1$ \\
			\hline
		\end{tabular}
	\end{center}
	Also, this discussion shows that in a general case, orienting covering space has maximum 4 components as follows 
	$P_0^0=(P_0,\mathcal O_{P_0}^0)$, $P_0^1=(P_0,\mathcal O_{P_0}^1)$, $P_1^0=(P_1,\mathcal O_{P_1}^0)$
	and
	$P_1^1=(P_1,\mathcal O_{P_1}^1)$,
	which respectively specify the orientations
	$\epsilon_0^0$, $\epsilon_0^1$, $\epsilon_1^0$ and $\epsilon_1^1$ on $M$.\\
	Let $\iota_a^b$ be an embedding of $P_a^b$ 
	into $P$ and let $U$ be an open neighborhood
	in $M$
	such that
	$\tilde{\pi}^{-1}(U)= V_1\cup V_2$
	and
	$\tilde{\pi}|_{V_i}$,
	$i=1,2$,
	be a homeomorphism.
	Then $\iota^*\circ \pi^*$ is an isomorphism from the connected component
	$(P_a, \mathcal{O}^b_{P_a})$
	to
	$ (M, \mathcal{O}_M)$.
	\subsection{Orientability of supermanifolds from point of view of vector bundles} \label{3.2}
	Let
	$M=(\tilde M,\mathcal O_M)$
	be a supermanifold and let $\mathcal J$ be the sub sheaf of ideals generated
	by odd elements. Therefore
	$\left( \dfrac{\mathcal J}{\mathcal J^2}\right) $ is a locally free
	$\mathcal{O/J}$- module.
	that
	$\mathcal J$
	is an ideal maximal ring
	$\mathcal O_M$.
	We show that if
	$M$
	be oriented, then its underlying manifold and vector bundle corresponding to
	$ \dfrac{\mathcal J}{\mathcal J^2}$
	both are oriented and vice versa.
	\begin{lemma}\label{3.1}
		Let
		$M^{m|n}=(\tilde M,\mathcal O_M)$
		be a supermanifold. Then there exists a vector bundle
		associated to $ \dfrac{\mathcal J}{\mathcal J^2}$.
	\end{lemma}
	\begin{proof}
		At any point
		$x$
		we consider
		the stalk
		$\mathcal{J}_x/\mathcal{J}^2_x$.
		Suppose
		$\xi_i$ for $1\leq i\leq n$, be odd coordinates in $U_{\alpha}$ a neighborhood of
		$x$.
		We consider vector space $E_x=\left\langle [\xi_i]_x+\mathcal J^2_x\right\rangle_\mathbb R $
		and put\\
		$E=\bigsqcup_{x \in \tilde M} E_x$. It can be shown
		$$\Xi=(E, \pi, \tilde M, \mathbb R^n)$$
		is a vector bundle such that its sheaf of sections is isomorphic with
		$\dfrac{\mathcal J}{\mathcal J^2}$.\\
		In the following, we will consider the local charts of a vector bundle
		$\Xi$.
		Set:
		\begin{align*}
			&h_\alpha: U_\alpha \times \mathbb R^n \to \pi^{-1}(U_\alpha)\\
			& (x,(v_1,...,v_n)) \mapsto \sum_{i=1}^{n}v_i([\xi_i]_x+\mathcal J^2_x)
		\end{align*}
		and
		\begin{align*}
			&h_\alpha^{-1}: \pi^{-1}(U_\alpha) \to U_\alpha \times \mathbb R^n \\
			& \sum_{i=1}^{n}\alpha_i([\xi_i]_x+\mathcal J^2_x) \mapsto (x,(\alpha_1,...,\alpha_n)).
		\end{align*}
		At any point
		$x \in U_\alpha \cap U_\beta$ we consider local coordinates
		$(x_1,...,x_m;\xi_1,...,\xi_n)$
		for
		$U_\alpha$
		and local coordinates $(y_1,...,y_m;\eta_1,...,\eta_n)$
		for
		$U_\beta$.
		So one has
		$\xi_i= \sum_{J\subset{1,...,n}}f_{iJ}\eta_J$.
		Therefore
		$$[\xi_i]_x +\mathcal{J}_x^2=\sum_{j=1,..., n}
		[f_{ij}]_x[\eta_j]_x + \mathcal{J}^2_x.$$
		The transition map is as follows:
		\begin{align*}
			h_{\beta}^{-1}\circ h_{\alpha} : & U_{\alpha\beta}\times \mathbb R^n \to U_{\alpha\beta}\times \mathbb R^n\\
			&(x,V)=(x,(v_1,...,v_n)) \mapsto \sum_{i} v_i([\xi_i]_x+\mathcal J^2_x)\\
			&\,\,\,\,\,\,\qquad\qquad\,\,\,\, \mapsto \sum_{i} v_i\left( (\sum_{j=1}^{n}	[f_{ij}]_x[\eta_j]_x)+\mathcal J^2_x\right)\\
			&\,\,\,\,\,\,\qquad\qquad\,\,\,\, =\sum_{j}\left(\sum_{i}  v_i 	[f_{ij}]_x\right) \left( [\eta_j]_x+\mathcal J^2_x\right)\\
			&\,\,\,\,\,\,\qquad\qquad\,\,\,\, \mapsto\left(x,\left(\sum_{i}  v_i	[f_{i1}]_x,...,\sum_{i}  v_i	[f_{in}]_x \right)  \right)\\
			&\,\,\,\,\,\,\qquad\qquad\,\,\,\, =\left( x,V(	[f_{ij}]_x)\right) 
		\end{align*}
		therefore
		$\Xi$ is a vector bundle.
	\end{proof}
	\begin{theorem}\label{4.1}
		The following propositions are equivalent: \\
		a) Supermanifold
		$M$ is orientable\cite{1}.\\
		b) The underlying manifold of
		$M$
		and vector bundle associated to
		$ \dfrac{\mathcal J}{\mathcal J^2}$
		are orientable.
	\end{theorem}
	\begin{proof}
		If
		$M$ is orientable then
		$\tilde M$ is clearly orientable. Regarding orientability of vector bundle $\dfrac{\mathcal J}{\mathcal J^2}$, since
		$M$
		is orientable, so its orienting covering space
		$P$
		has four connected components
		$P_i$,
		$i=1,2,3,4$, and oriented charts in each connected component induce an oriented trivializing atlas for vector bundle $E$. Thus the determinant of the transition map must be positive. According to Lemma
		\ref{3.1}
		transition map matrix is as follows;
		$$(f_{ij}(x))=\left( {\dfrac{\partial \xi}{\partial \eta}}\right)^\sim_x.$$
		Since
		$ \xi $
		and
		$\eta$ are odd coordinates of oriented charts in a connected component of $P$, therefore the determinant of the latter matrix is positive.
		Conversely, based on the assumption the vector bundle $E$ associated to
		$\dfrac{\mathcal J}{\mathcal J^2}$
		is orientable,
		so the
		supermanifold
		$( \tilde M ,\Gamma(\wedge( E )))$ is orientable.
		Since this supermanifold is isomorphic with $M$, thus $M$ is
		orientable too.
	\end{proof}
	\subsection{Semiorientable and nonorientable supermanifolds}
	Set
	$ P_a^b:= (a, b). P_0^0$ where (a, b) is an element of $\mathbb Z_2 \oplus \mathbb Z_2$.
	If these sets are distinct from each other, the covering space $P$ of the supermanifold
	$M$ has 4 components and
	$M$ is orientable. 
	But if for example
	$P_0^0=P_1^0$
	then
	$P_0^1=P_1^1$, in this case
	$P$
	has two components. When
	$P$
	has two connected components,
	$M$ is called
	\textit{semiorientable}.
	There exist three cases, based on whether the underlying manifold or associated vector bundle is orientable or not. If $P_0^0=
	P_a^b$ then this case is called 
	$(a,b)$-
	semi orientable.
	If
	$P$
	is connected then
	$M$
	is called nonorientable\cite{1}.\\
	Some examples of semi orientable and nonorientable supermanifolds are discussed in detail from \cite{1,5}.
	\section{Preimage orientation}
	In the case of supergeometry, similar to the classical geometry, if a supervector space 
	$V$
	is decomposed as 
	$V=U \oplus W$ then knowing the orientations of any two of these super vector spaces allows us to determine the orientation of the third super vector space as well. This orientation is called direct sum orientation. Now we introduce product orientation and preimage orientation.
	\begin{definition}[\textbf{product orientation}]
		Let
		$X^{m|n}$
		and
		$Y^{p|q}$
		be two oriented supermanifolds. Let
		$(x_1,...,x_m;\xi_1,...,\xi_n) $
		and
		$(y_1,...,y_p;\eta_1,...,\eta_q) $
		be coordinates represent the orientations on $X$
		and $Y$, respectively.
		Then the coordinates
		$$(x_1,...,x_m,y_1,...,y_p;\xi_1,...,\xi_n,\eta_1,...,\eta_q)$$ induces an orientation on
		$(X\times Y)^{m+p|n+q}$ which is called product orientation.
	\end{definition}
	\begin{definition} 
		Let
		$X$
		and
		$Y$
		be two oriented supermanifolds and $\psi: X \to Y$.
		At each point $x\in\tilde{X}$, in the given positively oriented bases for $T_xX$ and $T_yY$
		it is said that
		$\psi$
		is an orientation-preserving morphism whenever the matrix of the linear map
		$d\psi_x: T_xX \to T_yY$
		has a positive determinant at point $x$. Otherwise, it is called an orientation-reversing morphism.
	\end{definition}
	Let X and Y be supermanifolds and let Z be subsupermanifold of Y.
	Suppose the morphism $\psi:X\to Y$ is transversal to Z, abbreviated  $\psi \pitchfork Z$. In \cite{4} it is shown that
	$W=\psi^{-1}(Z)$ is subsupermanifold of X. 
	\begin{theorem} \label{4.3}
		If
		$X$, $Y$
		and
		$Z$ are oriented, then orientations on X, Y,  and Z induce an orientation on W.
	\end{theorem}
	\begin{proof}
		Since
		$\psi \pitchfork Z$,
		at any point
		$x \in \tilde W$
		we have
		\begin{equation}
			d\psi_xT_xX + T_yZ=T_yY
		\end{equation}
		where
		$y=\tilde{\psi}(x)$.
		Now let
		$W_x^c$
		be a complement vector space of
		$T_xW$
		in
		$T_xX$. Since
		$T_xW$
		is the preimage of
		$T_yZ$
		under the linear map
		$d\psi_x: T_xX \to T_zY$, we have
		\begin{equation}\label{6}
			d\psi_x W_x^c \oplus T_yZ=T_yY.
		\end{equation}
		Since Y and Z  are oriented, by \ref{6} one has a direct 
		sum orientation on $d\psi_x W^c_x$.
		$d\psi$ is monomorphism on $W^c$.
		Thus one gets an orientation on $W^c$.
		Now by the following equality:
		\begin{equation}
			W_x^c \oplus T_xW=T_xX
		\end{equation}
		one gets a direct sum orientation on $T_xW$. Let $(U_x, \epsilon_x)$ be a positively oriented chart on a neighborhood of $x$. Since $x$ is an arbitrary point, the neighbours $U_x$ with the specified orientations.\\
		constitute a connected component
		of oriented covering of W. 			In the end, for completing the proof, we show that W has a unique orientation.
		Let $(U_x,\epsilon_x)$ and  $(V_x, \epsilon'_x)$ are two intersected positively oriented charts on a neighborhood of $x$ in $X$. According to construction W \cite{4},
		since
		$U_x \cap V_x\neq 0$ 
		then by gluing conditions, the preimage orientation on the intersections $\epsilon_x$ and $\epsilon'_x$ must coincide.
	\end{proof}
	This orientation is called preimage orientation.
	\begin{proposition}
		Let
		$X$,
		$Y$
		and
		$Z$
		be oriented supermanifolds, and let Z be a sub supermanifold 
		of Y with $dimX+ 
		dimZ= dimY.$ In addition
		$\tilde X$
		is compact and
		$\tilde Z$
		is a closed submanifold of
		$\tilde Y$. If $\psi:X\to Y$ is a morphism transversal to Z, then it is possible to give an orientation pair of signs to the points of
		$W=\psi^{-1}(Z)$.
	\end{proposition}
	\begin{proof}
		Since $codim W= codimZ $ \cite{4}, by assumption $dimX+ dimZ= dimY$, we have $dim(W)=0|0$. So the complementary space of $W_x:= T_xW$, for $x\in W$, is
		$W_x^c$=$T_xX$. By the previous theorem we have the following relation:
		\begin{equation} \label{4.4}
			d\psi_xT_xX \oplus T_yZ=T_yY.
		\end{equation}
		According to the
		\begin{equation}
			d\psi_x^{-1} (d\psi_x T_xX) \oplus T_xW^{0|0}=T_xX
		\end{equation}
		By \ref{4.4} we have a direct sum orientation
		on $d\psi_x T_xX$. Since $d\psi_x$ is
		monomorphism, we get an orientation on $T_xX$.
		Denote this orientation by $\epsilon'_x$. By assumption 
		X has an orientation. It specifies an orientation
		on $T_xX$. Denote this orientation by $\epsilon_x$.
		Let $\epsilon_x=(\epsilon_0, \epsilon_1)$,
		where $\epsilon_0 (\epsilon_1)$ is the orientation induced by $\epsilon_x$ on even (odd) subspace of $T_xX$.
		Therefore $\epsilon'_i= \delta_{i, x} \epsilon_i$ for $i=0, 1$, where $\delta_{i, x}\in\{-1, +1\}.$
		Thus, the orientation sign at any point is one of the four ordered pairs
		$(\pm1, \pm1)$.
	\end{proof} 
	The pair $(\delta_{0, x}, \delta_{1, x})$ is called orientation pair of signs at $x$.
	\section{Oriented intersection pair}
	\subsection{Intersection pair}
	In this section, we will define the intersection number for transversal sub supermanifolds.
	\begin{definition}[\textbf{Oriented intersection pair}]\label{Def5.1}
		Let
		$X$
		be a compact supermanifold and
		$f: X \to Y$
		be a morphism between supermanifolds that transversal to a closed sub supermanifold
		$Z$
		in
		$Y.$ If
		$dim X+dim Z=dim Y$ then the underlying space
		$f^{-1}(\tilde Z)$
		is a finite set of points of
		$\tilde{X}$.
		The oriented intersection pair
		of $f$ and Z is denoted by $I(f, Z)$
		and is defined as follows:
		$$I(f,Z):=\sum_{x \in \tilde{f}^{-1}(Z)}(\delta_{0,x}, \delta_{1,x}).$$
		where $(\delta_{0,x}, \delta_{1,x})$ is orientation pair of signs at $x$, see the end of previous section.
	\end{definition}
	\begin{remark}
		Let X, Y, Z, and $f$ be as above and let $(U,{\epsilon_X})$ and $(V,{\epsilon_Y})$ be positively oriented local charts on X, Y respectively.
		We show that orientation pair of signs at any point
		$x \in \tilde W$, where $W=f^{-1}(Z)$, can be computed as follows:
	\end{remark}
	Let Jacobian matrices of coordinate representation of the map $f\times\iota: X\times Z \to Y\times Y$
	be as follows;
	$y=\tilde f(x) \in \tilde Z$.
	According to the
	\begin{equation}
		\begin{pmatrix}
			\begin{matrix} A & 0 \\ 0 & C \end{matrix} &\vline& 
			\begin{matrix} E & 0 \\ 0 & G \end{matrix} \\
			\hline
			\begin{matrix} B & 0 \\ 0 & D \end{matrix} &\vline& 
			\begin{matrix} F & 0 \\ 0 & H \end{matrix}
		\end{pmatrix}
	\end{equation}
	where $\iota$ is the inclusion map. Now consider the matrices
	$J_0= (\tilde{A}, \tilde{C})$ and 
	$J_1= (\tilde{F}, \tilde{H})$ where 
	$\tilde{A} ( or \tilde{C}, etc)$ is the matrix over 
	$\mathcal{O}_{X\times Z}/\mathcal{J}_{X\times Z}$  obtained by applying 
	the map $\mathcal{O}_{X\times Z}\to \mathcal{O}_{X\times Z}/\mathcal{J}_{X\times Z}$
	to each entry of A (or F, etc).
	It is seen that $\delta_{0, x}= sgn\,detJ_0(x, z)$ and $\delta_{1, x}= sgn\,detJ_1(x, z)$, where $z=\tilde{f}(x).$
	\begin{proposition} \label{5.3.0}
		Let X, and Y  be two supermanifolds and let $f: X\to Y$ be a morphism and Z
		a sub supermanifold of Y. If
		$f \pitchfork Z$
		then
		$\tilde {f} \pitchfork \tilde {Z}$
		and
		$I(\tilde{f}, \tilde{Z})= \sum_{x} \delta_{0,x}.$
	\end{proposition}
	\begin{proof}
		By theorem
		\ref{4.3}
		since
		$X$
		is oriented
		we take $\epsilon_0$
		to be the even-part orientation and
		$\delta_{0,x}$
		is obtained from the preimage orientation of the fixed orientation on
		$\tilde{X}$
		. We also know if
		$f$
		is transversal with
		$Z$
		then
		$\tilde f$
		is transversal with
		$\tilde Z$ \cite{4}. As a result, we have, by definition
		\begin{equation}
			d\tilde{f}_x T_x\tilde{X} \oplus T_y{\tilde{Z}}=T_y\tilde{Y}
		\end{equation}
		so the orientation on
		$T_y {\tilde{Y}}$
		and
		$T_y {\tilde{Z}}$
		induce a direct sum orientation on
		$d\tilde{f}_x T_x\tilde{X}$. And since
		$d\tilde{f}_x$
		is a monomorphism,
		then an orientation on
		$T_x {\tilde{X}}$
		corresponds to a preimage orientation, which if we denote this orientation by
		$\epsilon'_0$
		then
		$\epsilon'_0= \delta_{0,x}\epsilon_0 $
		, so
		$I(\tilde{f},\tilde{Z})=\sum_{x}\delta_{0,x}.$
	\end{proof} 
	To determine the second component of the intersection pair, we can use vector bundles  point of view given in the section
	\ref{3.2}. In what follows, we show that the second component of the intersection pair can be calculated as a classical term. To this end, let start with introducing some necessary notations and propositions and then we show that the second component of the intersection pair can be calculated by a homotopy invariant classical term. 
	\begin{definition}
		Let $M$, $N$ and $P$ be (smooth) manifolds and let $P$ be an embedded submanifold of $N$.  
		Assume $E_M$ and $E_N$ are two oriented vector bundles over $M$ and $N$ respectively. In addition, let $E_P$, the restriction of $E_N$ to $P$, be equipped with an orientation.
		\begin{enumerate}[label=(\alph*)] 
			\item \label{5.12}
			Let $\phi: E_M\to E_N$ be a bundle homomorphism
			which covers $\psi: M\to N$.  We say
			$\phi$ is transversal to $E_P$  if for each
			$x\in \psi^{-1}(P)$, one has
			$\phi(E_M)_x + (E_P)_{\psi(x)}= (E_N)_{\psi(x)}$.
			\item  \label{5.13.b} Let $\phi_x$ be injective. For $x\in\psi^{-1}(P)$, orientation sign of $\phi_x$ is denoted by $\delta_{\phi_x}$ and 
			equal to $\pm 1$.It is equal to +1 (res. -1) if the image of a positive base for $(E_M)_x$, under $\phi_x$, with a positive base of 
			$(E_P)_{\psi(x)}$ make a positive (res. negative) base for $(E_N)_{\psi(x)}$.	
			\item  Let $M$, $N$ and $P$ be smooth manifolds. In addition, assume $M$ is compact and P is a closed submanifold of $N$. Moreover $\psi$ is transversal to $P$ and 
			$dim M+ dim P= dim N$. Then $\psi^{-1}(P)$ is a finite set of points. The oriented intersection number of $\phi$ and $E_P$ is denoted by $I(\phi, E_P)$ and defined as follows: 
			$I(\phi, E_P)= \sum \delta_{\phi_x}$ where summation sign denotes a sum  over $x\in \psi^{-1}(P)$.
		\end{enumerate}
	\end{definition}
	Let $f: X \to Y$ be a morphism between supermanifolds. Then we have $f^*: \mathcal O_Y \to \mathcal O_X$, 	this morphism induces the following bundle morphism:
	$$\bar f^*:\Xi_{Y}\to \Xi_X$$
	where $\Xi_{Y}$ is vector bundle associated to $(\mathcal{J}/\mathcal{J}^2)_{Y}$, c.f. section \ref{3.2} On the dual bundles we have the following bundle morphism:
	$$\bar{f^*}^*: \left( {\dfrac{\mathcal J}{\mathcal J^2}}\right)^* _{X}\to \left( {\dfrac{\mathcal J}{\mathcal J^2}}\right)^*_{Y }.$$ 
	In the next proposition, we show that the dual sheaf of the quotient sheaf
	$\left( \dfrac{\mathcal J}{\mathcal J^2}\right) _X$
	at the point
	$x$
	is isomorphic  with the odd subspace of the tangent space
	$T_xX$
	. For this purpose, we introduce some useful morphisms. \\ 
	By $i_x$ we mean  the morphism 
	\begin{align}
		&i_x:\left( \left({\dfrac{\mathcal J}{\mathcal J^2}}\right) _X \right)_x ^* \to ({T_x}X)_1\notag\\
		&\qquad \qquad \qquad \quad \,\varphi \mapsto \tilde {\varphi}\label{15}
	\end{align}
	where the derivation $\tilde{\varphi}$ is defined as follows:
	\begin{equation}\label{16}
		\tilde{\varphi}(f):=\varphi((f-\tilde{f})+\mathcal J^2_X).
	\end{equation} Morphisms
	$i_z$
	and
	$i_y$
	are defined similarly.
	
	\begin{proposition}
		Morphisms
		$i_x$
		,
		$i_z$
		, and
		$i_y$
		defined above are isomorphisms.
	\end{proposition}
	\begin{proof}
		Every derivation
		$v \in ({T_x}X)_1$
		clearly defines a linear function
		$\varphi_v$
		on
		$\mathcal J_X$
		that is zero on
		$\mathcal J^2_X$.
		Conversely, for every
		$\varphi \in \left(\left(\dfrac{\mathcal J}{\mathcal J^2} \right)_{X}\right) _{x}^*$
		we have a tangent vector
		$v_\varphi$
		at
		$x$
		which is defined as
		$$v_\varphi(f)=\varphi\left( (f-\tilde f)+\mathcal J^2_X\right).$$
		The linearity condition for
		$v_\varphi$
		is clear. To check the Leibniz's property we have
		\begin{align*}
			v_\varphi(fg)=&\varphi\left( (fg-\tilde f\tilde g)+\mathcal J^2_X\right)\\
			=&\varphi\left( \left( (f-\tilde f)(g-\tilde g)+\tilde f(g -\tilde g)+(f-\tilde f)\tilde g\right) +\mathcal J^2_X\right)\\
			=&\tilde fv_\varphi(g)+(-1)^{p(f)p(g)}\tilde g v_\varphi(f).
		\end{align*}
		Therefore
		$v_\varphi$
		is a derivation. Also
		$v_{\varphi_v}$
		is
		equal to
		$v$.
	\end{proof}
	
	\begin{proposition} \label{5.9.0}
		Suppose
		$H:X \to Y$
		is a morphism between supermanifolds. Then for the bundle morphism
		\begin{align*}
			&(\bar{H^*})^*:\left(\left( \dfrac{\mathcal J}{\mathcal J^2} \right)_X\right)_x ^* \to \left(\left( \dfrac{\mathcal J}{\mathcal J^2} \right)_Y \right)_y ^*
		\end{align*}
		the following diagram is commutative.
		\begin{equation} \label{5.8.0}
			\begin{CD}
				({T_x}X)_1 @>d{H}>> ({T_y}Y)_1 \\
				@A{i_x}A{\rotatebox{90}{$\cong$}}A @A{\rotatebox{90}{$\cong$}}A{i_y}A \\
				\left( \left( \dfrac{\mathcal J}{\mathcal J^2} \right) _X \right)_x^* @>>(\bar{H^*})^*> \left( \left( \dfrac{\mathcal J}{\mathcal J^2} \right) _Y \right)_y^*
			\end{CD}
		\end{equation}
	\end{proposition}
	\begin{proof}
		For this purpose, we show
		$(d{H})_1 \circ i_x= i_y \circ (\bar{H^*})^*$.
		For
		$g \in \mathcal{O}_Y$
		on the left-hand side, we have:
		$$(d{H})_1 \circ i_x(\varphi)(g)=(d{H})_1 \circ {v_\varphi}(g)={v_\varphi} \circ H^*(g)={v_\varphi}(f) $$
		where
		$f=H^*(g) \in \mathcal{O}_X$.
		On the other hand, according to \ref{16} we have:
		\begin{align*}
			i_y \circ (\bar{H^*})^*(\varphi)(g)=&\left( \widetilde{(\bar{H^*})^*(\varphi)}\right) (g)\qquad\qquad\\
			=&\left( \varphi \circ (\overline{H_t^*})^*\right) \left( (g-\tilde{g})+{\mathcal J_Y}^2\right)\qquad\qquad\\ =&\varphi\left( {{{H}^*}}(g+{\mathcal J_Y}^2)-{{H}^*}(\tilde{g}+{\mathcal J_Y}^2)\right) \\
			=&\varphi\left( ({{H}^*}(g)-{{H}^*}(\tilde{g}))+{\mathcal J_Y}^2\right)\quad\quad \quad \\
			=&\varphi((f-\tilde f)+{\mathcal J_Y}^2)\\
			=&v_{\varphi}(f)
		\end{align*}
		The second equality holds
		according to \ref{5.8.0}
		and fourth equality holds
		by substituting $\widetilde{H^*}(g)$
		for $H^*(\tilde{g})$.
		therefore,  up to isomorphism, the morphisms
		$d{H}$
		and
		$(\bar{H^*})^*$
		are equivalent.
	\end{proof}
	\begin{proposition} \label{5.11.}
		Suppose
		$X$, $Y$, $Z$
		and f are as in Definition \ref{Def5.1}.
		Let $E_X$, $E_Y$ and $E_Z$ be canonical 
		vector bundle associated to X, Y, and Z
		respectively. Since f is transversal to Z 
		and by the proposition \ref{5.9.0}, it can be seen 
		that at each $x \in f^{-1}(Z)$, the image of a positive base for $(E_X)_x$, under ${\bar{f^*}}^*$, with a positive base for $(E_Z)_{f(x)}$ make a positive (res. negative) base for $(E_Y)_{f(x)}$ if $\delta_{1, x}>0$ (res. $\delta_{1,x}<0)$.
	\end{proposition}
	The next proposition provides a classical interpretation for oriented intersection number $I(f, Z)$, see Def. \ref{Def5.1}. Indeed, according to Def. \ref{5.13.b} and Prop.\ref{5.3.0} and prop. \ref{5.11.}, we have the following proposition.
	\begin{proposition} \label{5.14.}
		Let $X$, $Y$, $Z$ and f be as Def. \ref{Def5.1}. then $I(f, Z)= (I(\tilde{f}$, $\tilde{Z}), I({\bar{f^*}}^*, E_Z))$.
	\end{proposition}
	Based on the next subsection, the concept of transversality and oriented intersection pair can be generalized to bundle homomorphisms and subbundles.
	\subsection{Homotopy invariant}
	In this section, we discuss the issue of homotopy invariance of intersection pairs in supermanifolds.
	\begin{definition} \label{5.3.}
		Suppose
		$f: X^{m|n} \to Y^{p+m|q+n}$
		and
		$g:X^{m|n} \to Y^{p+m|q+n}$
		are morphisms between supermanifolds. We say
		$f$
		and
		$g$
		are homotopic, and denote it by
		$f \stackrel{H}{\sim} g$ if there exists a morphism
		$$H:X \times \mathbb{R}^{1|0} \to Y.$$
		One can also be viewed as a varying family of morphisms $\{H_t : X \to Y \,|\, t \in \mathbb R\}$,
		where for ${x \in \tilde X}$,
		$\tilde H_t (x) = \tilde H(x, t) $. Therefore
		$H \circ i_0=f$
		and
		$H \circ i_1=g$,
		where
		$ i_0\, \text{and}\,\, i_1: X^{m|n} \to X^{m|n} \times \mathbb{R}^{1|0}$
		are morphisms defined as follows:
		$$	\begin{tabular}{cc}
			$\begin{cases}
				\tilde i_0: \tilde X \to \tilde X \times \mathbb{R}^{1}\\
				\,\,\,\,\,\,\,\,\,\,\,\, x \mapsto (x,0)\\
			\end{cases}$
			&$\begin{cases}
				\tilde i_1: \tilde X \to \tilde X \times \mathbb{R}^{1}\\
				\,\,\,\,\,\,\,\,\,\,\,\, x \mapsto (x,1)\\
			\end{cases}$
		\end{tabular}$$
		also, the morphisms $i_0^*$ and $i_1^*$ are defined locally in terms of the coordinate system \\
		$(x_1, \ldots, x_m, t; e_1, \ldots, e_n)$, where $t$ is a coordinate to the supermanifold $\mathbb{R}^{1|0}$, as follows:
		\begin{center}
			\begin{tabular}{cc}
				$\begin{cases}
					i_0^*:f(x_1,\dots,x_m,t ) \mapsto f(x_1,\dots,x_m,0)\\
					i_0^*:e_j \mapsto e_j
				\end{cases}$&	$\begin{cases}
					i_1^*:f(x_1,\dots,x_m,t ) \mapsto f(x_1,\dots,x_m,1)\\
					i_1^*:e_j \mapsto e_j
				\end{cases}.$
			\end{tabular}
		\end{center}	
	\end{definition}
	\begin{theorem}[\textbf{Homotopy invariance of oriented intersection pair}] \label{5.4 Homotopy}
		Let
		$X^{m|n}$, $Y^{p+m|q+n}$ and
		$Z^{p|q}$ are supermanifolds.
		Consider two morphisms
		$f$
		and
		$g$, 
		from the supermanifold
		$X^{m|n}$
		to
		$Y^{p+m|q+n}$ such that both are transversal to a sub supermanifold
		$Z^{p|q}$
		of
		$Y$. Let
		$f$ and $g$,
		are homotopic and $H$ is a homotopy between them.
		If $H$ is transversal to $Z$, then
		$$I(f,Z)=I(g,Z).$$
	\end{theorem}
	\begin{proof}
		Consider $H:X \times \mathbb{R}^{1|0} \to Y$ as a homotopy 
		between $f$ and $g$ such that $f=H_{|_{X \times \{0\}}}$ and $g=H_{|_{X \times \{1\}}}$. For $t=0$
		$$\tilde H_0^{-1}(\tilde Z)=\tilde f^{-1}(\tilde Z)$$
		and
		for $t=1$, 
		$$\tilde H_1^{-1}(\tilde Z)=\tilde g^{-1}(\tilde Z)$$
		are two finite sets. In addition,
		$\tilde{H}^{-1}(\tilde{Z}) \cap \tilde{X} \times \tilde I$
		is a compact one-dimensional
		submanifold of 
		$\tilde{X}\times I$
		with boundary equal to 
		$\tilde{H}_0 \cup \tilde{H}_1$. 
		According to 
		Prop. \ref{5.3.0} and a classical result (see \cite{2}, Obsevation in Ch.3)
		we have $\sum \delta_{0, p_i}= \sum \delta_{0, q_j}$
		where the first sum is on the elements of $\tilde{H}_0$ and
		the second sum is on the elements of $\tilde{H}_1$.\\
		 Set 
		$h_t= H|_{X\times \{t\}}$. Then from the diagram \ref{5.8.0}  together with the classical reasoning
		as in the previous paragraph, it follows that the intersection orientation numbers of
		${\bar{h_0^*}}^*$ 
		and 
		${\bar{h_1^*}}^*$ are equal. So the second components of the intersection pair $f$ and $g$ are equal.
		Thus $I(f, Z)= I(g, Z).$
	\end{proof}
	\begin{remark}
		Let
		$X^{m|n}$
		, 
		$Y^{p+m|q+n}$ and
		$Z^{p|q}$ be supermanifolds.
		If $f:X \to Y$ such that $f$ is not transversal to $Z$.
		Then by genericity property of transversality \cite{4} there exists a morphism 
		$g$ which is homotopic to $f$ and $g \pitchfork Z$. Set
		$$I(f,Z):=I(g,Z).$$
		Theorem \ref{5.4 Homotopy} shows that the last definition is well-defined.
	\end{remark}
	In this section, we investigated the homotopy properties of the oriented intersection pair of supermanifolds.
	\subsection {Euler-Poincare characteristic pair}
	\begin{definition} \label{5.18.}
		Let
		$X=(\tilde X,\mathcal{O}_X)$
		be a supermanifold of dimension
		$m|n$.
		By diagonal of $X\times X$ we mean sub supermanifold
		$\Delta=(\tilde{\Delta},{\tilde{\iota}}_*\mathcal{O}_X)$ 
		where $\tilde{\Delta}$ is diagonal of 
		$\tilde{X} \times \tilde X$ and
		$\tilde{\iota}:\tilde{X} \to \tilde{\Delta}$ 
		is natural homeomorphism.
	\end{definition}
	Let $X$ and $Z$ be subsupermanifolds of $Y$. $X$ and $Z$ is said to be transversal if the inclusion map $j:X\to Y$ is transversal to $Z$. If $dim X + dim Z= dim Y$ then by $I(X, Z)$ we mean $I(j, Z)$.
	\begin{definition}[\textbf{Euler-Poincare characteristic pair}]\label{euler}
		Let $X$ be an oriented supermanifold 
		and let $\Delta$ be the diagonal of $X\times X$, see Def. \ref{5.18.}.
		The Euler- Poincaré characteristic pair of $X$ is denoted by 
		$\chi(X)$ and defined as follows
		$\chi(X)=I(\Delta, \Delta)$
	\end{definition}
	\begin{remark}
		It is worth noting that in common geometry
		$I(\Delta, \Delta)$ is equal to $L(id_X)$ the Lefschetz number of
		identity map $id: X\to X$. The Lefschetz theorem shows that this number is equal to the Euler- Poincare characteristic number of $X$.
	\end{remark}
	In the next section, we will investigate intersection number and Euler-Poincare characteristic pair in the category of
	$\Pi$-symmetric supermanifolds.
	This category and its morphism have been completely investigated in\cite{4}.
	
	\section{$\Pi$-symmetric supermanifolds}
	In this section, we will examine the intersection pair in the special case of $\Pi$-symmetric supermanifolds (see appendix \ref{A.3}). At the end of this section, we will calculate the Euler pair for an example of $\Pi$-symmetric Grassmannian supermanifolds.
	\subsection{Intersection pair of $\Pi$-symmetric supermanifold}
	At first, we discuss the intersection pair in the category of $\Pi$-symmetric supermanifolds.
	\begin{definition}
		Let $Y$
		is a
		$\Pi$-symmetric supermanifold, there exists one
		$\Pi$-symmetric morphism $P_Y:\mathcal T_Y \to \mathcal T_Y$ on the vector fields space of
		$Y$
		that induces at each $y\in\tilde{Y}$ a
		$\Pi$-symmetric morphism
		$P_Y^y: T_yY \to T_yY$ \cite{4}
		If Y is oriented, then $P_Y$ is 
		called orientation-preserving
		whenever for each y, $P^y_Y$ is
		orientation-preserving map.
	\end{definition}
	\begin{lemma} \label{2.3}
		Let X, Y,  and Z be objects in the category of $\Pi$-symmetric supermanifolds, see \cite{4}, such that $\tilde{X}$ is compact and $\tilde{Z}$ is a closed 
		sub supermanifold of $\tilde{Y}$ and let $f: X\to Y$ be a morphism in this category. If X, Y, and Z are oriented and $P_X$, 
		$P_Y$ and $P_Z$ are orientation-preserving
		morphisms, in addition
		$f \pitchfork Z$
		and
		$\dim X + \dim Z =\dim Y$
		then at any point
		$ x \in \tilde{\psi}^{-1}(\tilde Z)$ for the entries of orientation sign we have
		$\delta_{0,x}=\delta_{1,x}.$
	\end{lemma}
	\begin{proof}
		By
		$f\pitchfork Z$
		and
		$\dim X + \dim Z =\dim Y$
		at any point
		$x \in \tilde X$
		we have
		\begin{equation*}
			df_xT_xX \oplus T_yZ=T_yY.
		\end{equation*}
		Since
		$Y$
		is a
		$\Pi$-symmetric supermanifold, there exists one
		$\Pi$-symmetric morphism$P_Y: \mathcal T_Y ~\to~\mathcal T_Y$ on the vector fields space of
		$Y$
		that induces at each $y\in\tilde{Y}$ a
		$\Pi$-symmetric morphism
		$P_Y^y:T_yY\to T_yY$\cite{4}
		If Y is oriented, then $P_Y$ is 
		called orientation preserving
		whenever for each $y$, $P^y_Y$ is
		orientation preserving map.
		If $(y_1, ..., y_{m+ p})$ is a slice coordinates for $\tilde{Z}$ in 
		$\tilde{Y}$ around a point $y\in \tilde{Z}$,
		then
		$$T_yY=\left\langle \dfrac{\partial}{\partial y_1},...,\dfrac{\partial}{\partial y_{m+p}}; P_Y^y\dfrac{\partial} {\partial y_1},...,P_Y^y\dfrac{\partial}{\partial y_{m+p}}\right\rangle $$
		and
		$$T_yZ=\left\langle \dfrac{\partial}{\partial y_1},...,\dfrac{\partial}{\partial y_{p}}; P_Y^y\dfrac{\partial}{\partial y_1},...,P_Y^y\dfrac{\partial}{\partial y_{p}}\right\rangle. $$
		Since
		$X$
		is a
		$\Pi$-symmetric supermanifold, for each $x\in\tilde{f}^{-1}(y)$ we have
		$$T_xX=\left\langle \dfrac{\partial}{\partial x_1},...,\dfrac{\partial}{\partial x_{m}}, P_X^x\dfrac{\partial}{\partial x_1 },...,P_X^x\dfrac{\partial}{\partial x_{m}}\right\rangle .$$
		The map f is $\Pi$-symmetric morphism. Thus the following diagram commutes:
		\begin{equation}
			\begin{CD}
				T_xX@>df_x>>T_{\tilde{f}(x)}Y@.\qquad \\
				@ VP_XVV @ VV P_Y V @.\\
				T_xX@>>df_x>T_{\tilde{f}(x)}Y@.
			\end{CD}.
		\end{equation}
		Thus the matrix of the linear map
		$df_x$, in block form, is as follows:
		\begin{equation} \label{5.4}
			\begin{bmatrix}
				\alpha_{ij}& 0\\
				0 & \alpha_{ij}
			\end{bmatrix}
		\end{equation}
		where $df_x(\partial/\partial x_i)= \alpha_{ij}(\partial/\partial y_j).$
		Then according to the definition of orientation sign, \ref{5.4} implies that $\delta_{0,x}=\delta_{1,x}$. This completes the proof.
	\end{proof}
	\begin{remark}
		According to the lemma \ref{2.3}, in the special case of $\Pi$-symmetric supermanifolds, the intersection pair is as follows:
		$$I(f,z)=(I(\tilde f, \tilde Z),I(\tilde f, \tilde Z)).$$
	\end{remark}
	\begin{proposition}
		According to the definition \ref{euler} for every $\Pi$-symmetric supermanifold $X$ we have 
		$$\chi(X)=(I(\tilde \Delta, \tilde \Delta),I(\tilde \Delta, \tilde \Delta)).$$
	\end{proposition}	
	
	In the next section, we present an example of $\Pi$-symmetric Grassmannian and calculate its Euler characteristic pair.
	\section{Grassmannian supermanifolds}
	In this section, we first take a brief look at the definition and construction of Grassmannian supermanifolds. Then, $\Pi$-Grassmannians are introduced.
	Finally, we compute the Euler- Poincaré characteristic pair for a specific case.\\
	The set of all complex sub supervector spaces of dimension
	$k|l$
	of $\mathbb{C}^{m|n}$
	is called supergrassmannian
	and denoted by
	$Gr(k|l,m|n)$. This is a supermanifold 
	that can be constructed by gluing some copies of superdomains as follows:
	Let
	$I=I_0|I_1$
	be two sets  such that
	$I_0 \subset \{1,\dots,m\}$
	and
	$I_1 \subset \{1,\dots,n\}$
	with
	$|I_0|=k$
	and
	$|I_1|=l$. Consider the standard superdomains
	$U_I$
	with coordinates
	$(x^I,\xi^I)$, and
	$$U_I={\mathbb{C}^{k(m-k)+l(n-l)|k(n-l)+l(m-k)}}.$$
	Each $U_I$ is labeled by a matrix, say $Z_I$, in a standard format which the columns with indices in I together form a unit matrix as follows:
	\[
	Z_I=\begin{array}{|c|c|c|c|}	
		\hline
		I_k & (x^I)_{k \times (m-k)} & \mathbf{0}_{k \times l} & (\xi^I)_{k \times (n-l)} \\ \hline
		\mathbf{0}_{l \times k} & (\xi^I)_{l \times (m-k)} & I_l & (x^I)_{l \times (n-l)} \\ \hline
	\end{array}.
	\]
	Supergrassmannians are obtained by gluing together the latter superdomains. To introduce the gluing morphism, set
	$$Z_I={B_{IJ}}^{-1} Z_J$$
	where the matrix
	${B_{IJ}}$
	is obtained by selecting the columns in $Z_J$ with indices in I. Now a morphism that maps
	each $x^I$ or $\xi^J$ to the corresponding 
	entry in the last equation is the  
	morphism which glues superdomains 
	$U_I$ and $U_J$.
	\subsection{$\Pi$-symmetric Grassmannians}
	A $\Pi$-symmetric Grassmannian 
	$Gr^{\Pi}(k|k, m|m)$
	is a supermanifold that is constructed
	by gluing super domains
	$U_I= \mathbb{C}^{m(m-k)|m(m-k)}$
	with local coordinates
	$(x^I,\xi^I)$
	which $U_I$ is labeled by the following 
	matrix
	\[
	Z^{\Pi}_I=\begin{array}{|c|c|c|c|}
		\hline
		A & I_k & B & \mathbf{0} \\ \hline
		C & \mathbf{0} & D & I_l \\ \hline
	\end{array}.
	\]
	with $A=D=(x^I)_{k \times (m-k)}$ and $B= -C=(\xi^I)_{k \times (m-k)}$ (see \cite{9,10}).In the special case if
	$k=1$
	,
	$Gr^{\Pi}(1|1,m|m)$
	has complex projective space $\mathbb{C}P^{m-1}$ as its reduced manifold.
	\begin{example}
		In this example, we show that
		the Euler characteristic pair of 
		$Gr^{\Pi}(1|1, m|m)$ is equal to
		$(m, m)$.
		Since $Gr^{\Pi}(1|1, m|m)$ is a
		$\Pi$-symmetric supermanifold, see 
		the sentence just before this example,
		by lemma \ref{2.3} the two components of Euler characteristic pair of $Gr^{\Pi}(1|1, m|m)$ are equal. According to proposition \ref{5.3.0}
		the first component of this pair is equal to $\chi(\mathbb {CP}^{m-1})$; Euler-Poincare characteristic of the $\mathbb{CP}^{m-1}$ . Since $\chi(\mathbb{C}P^{m-1})=m$, our assertion holds.
	\end{example}
	\section{Conclusion}
	In this article, we first defined the intersection pair, using the concept of orientation in supermanifolds. This pair is obtained from the intersection number of the underlying manifolds and the intersection number of corresponding vector bundles. This interpretation reduces the calculation of the intersection number of the odd part to a calculation in classical geometry. In the end, we defined the Euler-Poincare characteristic pair in the category of $\Pi$-symmetric supermanifolds.
	\begin{appendix}
		\section{Appendix}
		\subsection{\,Some preliminaries of sheave and ringed spaces}
		In this section, some necessary concepts and results are introduced as follows:\\
		\textbf{Presheaf}.
		If $X$ is a topological space, then a presheaf is a correspondence, say $R$, that assigns to each open subset of $X$, say $U$, an algebraic structure $R(U)$. The correspondence $R$ must satisfy the following condition: if $V$ is an open subset of $X$ and $U$ is an open subset of $V$, then there exists a homomorphism $$r_{UV}: R(V) \to R(U)$$
		This homomorphism is called the restriction chart.
		Another important condition is that if $U \subset V \subset W$ are open subsets, then we have
		$$r_{UW}=r_{UV} \circ r_{VW}$$
		where $\circ$ denotes the composition of homomorphisms.\\
		\textbf{Sheaf}.
		A sheaf is a presheaf say $R$ with the following condition: If $U$ is an open subset and $\{U_{\alpha}\}_\alpha$ is an open cover of $U$ and $\{f_{\alpha}\}$ is  a family of elements 	$f_\alpha \in R(U_\alpha)$ such that for any $ \alpha,\beta$ with nonempty intersection $U_\alpha \cap U_\beta$ one has
		$$f_\alpha|_{U_\alpha \cap U_\beta} = f_\beta|_{ U_\alpha \cap U_\beta}$$
		then there exists a unique element $f \in R(U)$ such that $f|_{U\alpha}= f_\alpha$.\\
		\textbf{Ringed space}.
		Let $X$ be a topological space and $R$ be a sheaf of rings i.e. for every open subset $U$ of $X$, $R(U)$ is a commutative ring. Then $(X, R)$ is called a ringed space.  The elements of $R(U)$ are called the \textit{sections} of $R$ over $U$. We define the \textit{stalk} $R_x$ of $R$ at $x\in X$ to be the direct limit of the rings $R(U)$ for all open sets $U$ containing $x$ via the restriction charts. A ringed space is said to be a \textit{space} if the stalks are local, meaning that each stalk $R_x$ has a unique maximal ideal denoted by $\mathit{m}_x$.
		\\
		\textbf{Products}.
		The category of supermanifolds admits products. Let $X_{i}\, (1 \leq i\leq n)$ be spaces in the category. A ringed space $X$ together with projection charts $P_{i}$ : $X\rightarrow X_{i}$ is called a product of the $X_{i},$
		$$
		X=X_{1}\times\cdots\times X_{n},
		$$
		if the following is satisfied: for any ringed space $\mathrm{Y}$, the chart
		\begin{equation}\label{2.1..}
			f\mapsto(P_{1}\circ f,\ \ldots,\ P_{n}\circ f)
		\end{equation}
		from $\mathrm{H}\mathrm{o}\mathrm{m}(\mathrm{Y},\ X)$ to $\prod_{i}\mathrm{H}\mathrm{o}\mathrm{m}(Y,\ X_{i})$ is a bijection. In other words, every morphism $f$ from $\mathrm{Y}$ to $X$ are identified with $n$-tuples $(f_{1},\ \ldots,\, f_{n})$ of morphisms $f_{i}(\mathrm{Y}\rightarrow X_{i})$. This product is unique up to isomorphism .\\
		\textbf{Superdomain}.
		Let $U$ be an open set in $\mathbb{R}^p$, and let $\mathcal{C}^{\infty \, p|q}$ be the sheaf of smooth supercommutative $\mathbb Z_2$-graded rings on $U$ that assigns to each open $V \subseteq U$ the ring	$ C^\infty (V) \,\, [ \theta^1, ..., \theta^q ]$ i.e. the ring of polynomials in the variables $ \theta^1 ,..., \theta^q$, satisfying the relations $ \theta^i \theta^j= - \theta^j \theta^i$ and $\theta^{i^2}=0$ and coefficients in $C^{\infty}(V)$ the ring of smooth real valued functions on $V$. Then the pair $(U, \mathcal{C}^{\infty \, p|q} )$ is called a superdomain, denoted by $U^{p|q}$.
		\begin{remark} Each element of $\mathcal{C}^{\infty p|q}(V)$ can be written as 
			$
			\sum\limits_{I \subset\{1,\dots,q\}} f_I\theta^I,
			$
			where the $f_I \in C^{\infty}(V)$ and $\theta^{I}=\theta^{i_1} \cdots \theta^{i_r}$ if $I=\{i_1,\dots , i_r\}$ such that $\ i_1< \dots < i_r$.
		\end{remark}
		\subsection{Examples of semi and nonorientable supermanifolds}
		\begin{example}
			Let
			$\overset{..}{M}$
			be a Möbius strip. It can be considered as a nonorientable line bundle over a circle discussed in detail from \cite{5}
			in the appendix. Let denote the sheaf of sections of an exterior bundle of M by $\Gamma(\wedge M)$. Then the supermanifold
			$$ N^{1|1}=(S^1,\Gamma(\wedge \overset{..}{M}))$$
			is an example of (1,0)-semi-orientable supermanifold, because its underlying manifold is the orientable circle, and the line bundle over it is not orientable. Indeed,
			$\overset{..}{M}=\dfrac{I\times\mathbb R}{\sim}$, where $I= [0, 1]$ and $\sim$
			is an equivalence relation from
			under which any pairs $(0, v)$ and $(1, -v)$ are equivalent. Thus they specify points in quotient space on which orientation classes induced by coordinates
			$(x, t)$
			and
			$(x, -t)$
			are equal to each other.
			So $P^0_0$ is equal to $P^1_0$
			and therefore the orienting covering space of
			$N^{1|1}$
			has two components.
		\end{example}
		\begin{example}
			Let M be a Mobius strip. It is a two-dimensional
			nonorientable manifold. We consider $\overset{..}{M}\times \mathbb R$
			as a trivial vector bundle on M. Therefore
			$$K^{2|1}=(\overset{..}{M},\Gamma(\wedge(\overset{..}{M} \times\mathbb R)))$$
			is an example of a 
			(0,1)-semi orientable supermanifold. In any local chart, the following orientations
			are equal:
			\begin{align*}
				&[(x_1, x_2; \xi)]=[(x_1, -x_2; \xi)]\\
				& [(x_1, x_2; -\xi)]= [(x_1, -x_2; -\xi)].
			\end{align*}
			Thus the orienting covering space for
			$K^{2|1}$
			has two components.
		\end{example}
		\begin{example}
			In this example, we consider the Möbius strip $\overset{..}{M}$ as a line bundle over $S^1$. Let
			$\overset{..}{M}\oplus\overset{..}{M} $
			be the Whitney sum of two copies of $\overset{..}{M}$. This space can be considered as a linear vector bundle on the Möbius strip.
			Then the
			$$S^{2|1}=(\overset{..}{M},\Gamma(\wedge(\overset{..}{M}\oplus\overset{..}{M}))) $$
			is an example of (1,1) semi-orientable supermanifolds.\\
			Note that $\overset{..}{M}\oplus \overset{..}{M}$
			can be considered as the quotient space of
			$\mathbb{R}^3 $
			under the following equivalence relation:
			$$(x, y, z)\sim (x+n, (-1)^n y, (-1)^n z).$$
			Note that in any local chart orientation induced by the following coordinate systems are equal:
			$$ (x, y, z)  \text{
				,} (x, -y, -z).$$
			Therefore, the orienting covering space of
			$S^{2|1}$, say $P$, has two components 
			and $P_0^0= P_1^1$.
		\end{example}
		\begin{example}
			In this example, we introduce a nonorientable supermanifold. Let consider
			$$C^{3|2}=(S^1\times\overset{..}{M},\Gamma(\wedge(\overset{..}{M}\times\overset{..}{ M}\times \mathbb R)) )$$
			
			$N^{1|1}\times K^{2|1}$ as a product of two supermanifolds $N^{1|1}$ and $K^{2|1}$. If
			$(x_1, x_2, x_3; \xi_1, \xi_2)$
			be the local coordinate system for this supermanifold, then following orientation classes are equal
			$$[(x_1, x_2, x_3; \xi_1, \xi_2)]=[(x_1, x_2, -x_3; \xi_1, \xi_2)]= [(x_1, x_2, x_3;- \xi_1, \xi_2 )]= [(x_1, x_2,- x_3; -\xi_1, \xi_2)].$$
			So the orienting covering space is connected.
		\end{example}
		\subsection{$\Pi$-symmetric supermanifolds}\label{A.3}
		A $\Pi$-symmetric free
		$A$-module $S$ of rank $m|m$, is defined to be an $A$-module isomorphic to
		$$ A^m \oplus (\Pi A)^m.$$
		Now we consider
		$A$-module isomorphism
		$ P_S: S \to \Pi S$ such that
		$P_S^2=id_S$
		where
		$\Pi$
		be a parity change functor.
		We call
		$P_S$, a
		$\Pi$ -symmetry on
		$S$.
		Objects in the category of $\Pi$-symmetric supermanifolds are supermanifolds such that there exists a $\Pi$-symmetry on their tangent sheaf. From now, for convenience, by $\Pi$-symmetry on a supermanifold we mean a $\Pi$- symmetry on its tangent sheaf.
		Let $M$ and $N$ be $\Pi$- symmetric supermnifolds. Then the set of morphisms $Hom(M, N)$  consists of all morphism $\psi: M\to N$ such that the following diagram commutes
		\begin{equation}
			\begin{CD}
				\psi_*(\mathcal T_M)@> {\bar D(\psi)}>>\mathcal T_{N,M}@.\qquad \\
				@ V{P_M}VV @ VV{\bar P_N} V @.\\
				\Pi\psi_*(\mathcal T_M)@>>{\bar D(\psi)^\Pi}>\Pi\mathcal T_{N,M}@.
			\end{CD}
		\end{equation}
		where $P_M$ is $\Pi$-symmetry on $\mathcal{T}_M$ and $\bar{P}_N$ is a  morphism induced by $\Pi$-symmetry $P_N$ on $\mathcal{T}_{N,M}:= Der(\mathcal{O}_N \to \mathcal{O}_M)$ such that $\bar{P}_N(f\partial/\partial y):=fP_N(\partial/\partial y)$ for each $f \in \mathcal{O}_M$ and 
		$\partial/\partial y$ is a coordinate vetor field on $N$.
		
		In addition $\bar{D}(\psi): \psi_*(\mathcal{T}_M)\to \mathcal{T}_{N,M}$ is defined as follows:
		$\bar{D}(\psi){X}(g)= X(\psi^*(g))$ where $X \in Der(\mathcal O_M)$ and $g \in \mathcal O_N$.\\
		
		To prove the commutative of the diagram, we must act as follows;\\
		Let
		$(x_1,...,x_n;\theta_1,...,\theta_n)$
		and
		$(y_1,...,y_m;\xi_1,...,\xi_m)$
		are local coordinates 
		for $U_\alpha \subseteq M$ and $U_\beta \subseteq N$, respectively.
		$$(\bar{P}_N \circ \bar{D}(\psi)){X}(g)=({\bar D(\psi)^\Pi} \circ P_M) {X}(g).$$
		\begin{example}
			Let $X={\dfrac{\partial}{\partial x_i}}$ and $\psi^*(f)=g;$
			\begin{align*}
				&(\bar{P}_N \circ \bar{D}(\psi)){X}(f)=\bar{P}_N \circ (\bar{D}(\psi){\dfrac{\partial}{\partial x_i}}(f))=\bar{P}_N \circ ({\dfrac{\partial}{\partial x_i}})(g)=\\
				&(\bar{P}_N \circ ({\dfrac{\partial}{\partial x_i}}))(g)=	\bar{P}_N \circ (g_i {\dfrac{\partial}{\partial y_i}})(g)=g_i P_N({\dfrac{\partial}{\partial y_i}})(g)=g_i ({\dfrac{\partial}{\partial \xi_i}})(g)={\dfrac{\partial g }{\partial \theta_i}}
			\end{align*}
			On the other hand;
			\begin{align*}
				&({\bar D(\psi)^\Pi} \circ P_M) {X}(f)={\bar D(\psi)^\Pi} \circ (P_M(X)(f))={\bar D(\psi)^\Pi} \circ (P_M({\dfrac{\partial}{\partial x_i}})(f))=\\
				&\qquad \qquad \qquad {\bar D(\psi)^\Pi} \circ ({\dfrac{\partial}{\partial \theta_i}})(f))={\dfrac{\partial}{\partial \theta_i}}(\psi^*f)={\dfrac{\partial g}{\partial \theta_i}}
			\end{align*}
			In the following, we will consider $\phi:Y^{l|l} \to \mathbb{R}^{n|n}$ is $\Pi$-symmetric or the diagram below in a special case is commutative;
			\begin{equation}
				\begin{CD}
					\phi_*(\mathcal T_{{Y^{l|l}}})@> {\bar D(\phi)}>>\mathcal T_{{\mathbb{R}^{n|n}},{{Y^{l|l}}}}@.\qquad \\
					@ V{P_{Y^{l|l}}}VV @ VV{\bar P_{\mathbb{R}^{n|n}}} V @.\\
					\Pi\phi_*(\mathcal T_{_{Y^{l|l}}})@>>{\bar D(\phi)^\Pi}>\Pi\mathcal T_{{\mathbb{R}^{n|n}},{{Y^{l|l}}}}@.
				\end{CD}
			\end{equation}
			According to embedding $\phi$, there exists $U_\alpha \subseteq \mathbb{R}^{n|n}$ and $V_\beta \subseteq \mathbb{Y}^{l|l}$ such that;
			$$\phi_*:\mathcal{O}_{\mathbb{R}^{n|n}}({U_\alpha}) \to \mathcal{O}_{Y^{l|l}}({V_\beta})$$
			$$(x_1,\dots,x_n;\theta_1,\dots,\theta_n)\mapsto(x_1,\dots,x_l,0,\dots,0;\theta_1,\dots,\theta_l,0,\dots,0)$$
			We must prove;
			$$(\bar{P}_{\mathbb{R}^{n|n}} \circ \bar{D}(\phi)){X}(g)=({\bar D(\phi)^\Pi} \circ P_{Y^{l|l}}) {X}(g)$$
			For this purpose, let $X={\dfrac{\partial}{\partial x_i}}$ and $\phi^*(g)=f;$
			\begin{align*}
				&(\bar{P}_{\mathbb{R}^{n|n}} \circ \bar{D}(\phi)){X}(g)=\bar{P}_{\mathbb{R}^{n|n}} \circ (\bar{D}(\phi)({\dfrac{\partial}{\partial x_i}})(g))=\bar{P}_{\mathbb{R}^{n|n}} \circ ({\dfrac{\partial}{\partial x_i}}(\phi^*(g)))\\
				&\qquad \qquad \qquad =\bar{P}_{\mathbb{R}^{n|n}} \circ ({\dfrac{\partial f}{\partial x_i}})
				=f_i \bar{P}_{\mathbb{R}^{n|n}} ({\dfrac{\partial}{\partial \theta_i}})=f_i {\dfrac{\partial}{\partial \theta_i}}={\dfrac{\partial f}{\partial \theta_i}}
			\end{align*}
			On the other hand;
			$$({\bar D(\phi)^\Pi} \circ P_{Y^{l|l}}) {X}(g)={\bar D(\phi)^\Pi} \circ (P_{Y^{l|l}}({\dfrac{\partial}{\partial x_i}})(g))={\bar D(\phi)^\Pi}({\dfrac{\partial}{\partial \theta_i}})(g)$$
			$$
			={\dfrac{\partial}{\partial \theta_i}}(\phi^*(g))={\dfrac{\partial f}{\partial \theta_i}}
			$$
		\end{example}
	\end{appendix}
	\bibliographystyle{amsplain}
		
\end{document}